\documentclass{article}
\usepackage{amsmath}
\usepackage[margin=10truemm]{geometry}
\title{\huge Diophantine Equation $n(x^4 + y^4) = z^4 + w^4$}
\author{Seiji Tomita and  Oliver Couto}
\date{}
\begin{document}
\maketitle
\begin{abstract}
In this paper, we proved that there are infinitely many  parametric solutions of $n(x^4 + y^4) = z^4 + w^4.$\\
\end{abstract}

\centerline{\huge1. Introduction}
\vskip\baselineskip
In 2016 Izadi and Nabardi\cite{b} showed $x^4+y^4=2(u^4+v^4)$ has infinitely many integer solutions.\\
They used  a specifc congruent number elliptic curve namely $y^2 = x^3 - 36x^2$.\\
In 2019 Janfada and  Nabardi\cite{c} showed that a necessary condition for $n$ to have an
integral solution for the equation $x^4+y^4=n(u^4+v^4)$ and gave a parametric solution.\\
They gave the numeric solutions for $n=41,136,313,1028,1201,3281,\cdots$.

In 2020 Ajai Choudhry\cite{a}, Iliya Bluskov and Alexander James showed $(x_1^4+x_2^4)(y_1^4+y_2^4)=z_1^4+z_2^4$ has 
infinitely many parametric solutions.
They gave the numeric solutions for $n=17,257,626,641,706,1921,\cdots$.

Inspired by Choudhry's article, we showed if $n=(a^4+b^4)/2$ then $n(x^4 + y^4) = z^4 + w^4$ has infinitely many parametric solutions.\\
We consider n be an integer.
Using four different methods we showed $n(x^4 + y^4) = z^4 + w^4$ has infinitely many integer solutions.\\
First method is similar as Choudhry's one, second method is using another identity.\\
Moreover, we showed the third method through an example $17(x^4+y^4)=z^4+w^4$.\\
In addition, we showed the fourth method through an example $97(x^4+y^4)=z^4+w^4$.\\
Finally, we gave the numeric solutions for $n=2,8,17,41,82,97,113,136,137,146,178,193,\cdots$ with $n<1000$.

\newpage

\centerline{\huge2. First method}
\vskip\baselineskip

\textbf{Theorem 1}
$(x^4+y^4)(a^2+b^2)=z^4+w^4$ has infinitely many parametric solutions,\\
where $(a,b)=(\frac{m^2-n^2}{2},\frac{m^2+n^2}{2})$, $m,n$ are arbitrary integer and have same parity.
\vskip\baselineskip

Proof.

We use Brahmagupta's identity $$(A^2+B^2)(C^2+D^2)=(AC+BD)^2+(AD-BC)^2$$.

We obtain $$(x^4+y^4)(a^2+b^2) = (x^2a+y^2b)^2 + (x^2b-y^2a)^2$$

Hence we have to find the rational solution of $x^2a+y^2b =u^2$ and $x^2b-y^2a = v^2$.

Parametric solution of first equation is $(x,y)=(k^2-2mk+4a-3m^2, -k^2-2mk+4a-m^2$) where $a+b = m^2$. \\
Substitute it to second equation, we obtain

$$v^2 = (m^2-2a)k^4-4k^3m^3+(4am^2-2m^4)k^2+(-32am^3+32a^2m+12m^5)k+48a^2m^2-32a^3-34am^4+9m^6$$

Let $a = \frac{m^2-n^2}{2}$ and $U=\frac{1}{k}$, then we obtain
\begin{equation}
V^2 = (5m^4n^2+4n^6)U^4+(4m^5+8mn^4)U^3-2m^2n^2U^2-4m^3U+n^2
\end{equation}

Since equation $(1)$ has a point $Q(U,V)=(0,n)$,
then equation $(1)$ is birationally equivalent to an elliptic curve $E$.
\begin{align*}
E: Y^2-4m^3YX/n+(8nm^5+16mn^5)Y &= X^3-2m^2(n^4+2m^4)X^2/(n^2)+(-20m^4n^4-16n^8)X  \\
                                &+104m^6n^6+32m^2n^{10}+80m^{10}n^2
\end{align*}
                                     
Transformation is given, 
\begin{align*}
U &= \frac{2n^2X-4m^2n^4-8m^6}{Yn} \\
V &= \frac{n^4X^3-6n^6X^2m^2-12X^2m^6n^2+28n^8m^4X+32n^4m^8X+32m^{12}X-16m^9Yn+16n^9mY+16n^{12}X-104n^{10}m^6-32n^{14}m^2-80n^6m^{10}}{Y^2n^3} \\
X &= \frac{2nV+2n^2-4m^3U}{U^2} \\
Y &= \frac{4n^3V+4n^4-8n^2m^3U-4n^4m^2U^2-8m^6U^2}{U^3n}
\end{align*}

The point corresponding to point $Q$ is $P(X,Y)=( \frac{2m^2(n^4+2m^4)}{n^2}, \frac{-16m(-m^8+n^8)}{n^3})$. \\
According to the Nagell-Luts theorem, since the point $P$ is a point of infinite order then we can obtain infinitely many rational points on $E$.
Thus we can obtain infinitely many rational solutions of $(1)$. \\
Therefore $(x^4+y^4)(a^2+b^2)=z^4+w^4$ has infinitely many parametric solutions. \\
The proof is completed. \\

When $m,n$ have opposite parity, problem links to $h(x^4+y^4) = 2(z^4+w^4)$ where $a^2+b^2=h/2$. 
\vskip\baselineskip

Example

$(a,b)=(\frac{m^2-n^2}{2},\frac{m^2+n^2}{2})$ \\
$m,n$ are arbitrary integer and have same parity.

$2Q(U) = \frac{-2n^2m}{n^4+m^4}$
\begin{align*}
&x = -n^4+4m^2n^2+m^4 \\
&y =  n^4+4m^2n^2-m^4 \\
&z = (3n^4+m^4)m \\
&w = (n^4+3m^4)n
\end{align*}
\vskip\baselineskip

$3Q(U) = \frac{2mn^2(n^4+3m^4)}{n^8-8m^4n^4-m^8}$
\begin{align*}
&x = m^{12}+12n^2m^{10}-19m^8n^4+40n^6m^6+19m^4n^8+12n^{10}m^2-n^{12} \\
&y = -m^{12}+12n^2m^{10}+19m^8n^4+40n^6m^6-19m^4n^8+12n^{10}m^2+n^{12} \\
&z = m(m^{12}+41m^8n^4+27m^4n^8-5n^{12}) \\
&w = n(-n^{12}-41m^4n^8-27m^8n^4+5m^{12})
\end{align*}
\vskip\baselineskip

Numerical example: $(m,n)<10$ \\
\begin{table}[hbtp]
\centering
  \caption{Solutions of $(x^4+y^4)(a^2+b^2) = z^4+w^4$}
  \label{table:data_type}
  \begin{tabular}{cclllll}
    \hline
   m & n & $a^2+b^2$ & x  & y  &  z  & w   \\
    \hline \hline
    
3 & 1 & 41 & 29 & 11 & 63 & 61 \\
4 & 2 & 136 & 31 & 1 & 76 & 98 \\
5 & 1 & 313 & 181 & 131 & 785 & 469 \\
5 & 3 & 353 & 361 & 89 & 1085 & 1467 \\
6 & 4 & 776 & 209 & 79 & 774 & 1036 \\
7 & 1 & 1201 & 649 & 551 & 4207 & 1801 \\
7 & 3 & 1241 & 1021 & 139 & 4627 & 5463 \\
7 & 5 & 1513 & 1669 & 781 & 7483 & 9785 \\
8 & 2 & 2056 & 319 & 191 & 2072 & 1538 \\
8 & 6 & 2696 & 751 & 401 & 3992 & 5094 \\
9 & 1 & 3281 & 1721 & 1559 & 14769 & 4921 \\
9 & 5 & 3593 & 3509 & 541 & 18981 & 25385 \\
9 & 7 & 4481 & 5009 & 2929 & 30969 & 38647 \\
    \hline
  \end{tabular}
\end{table}

\newpage

\centerline{\huge3. Second method}
\vskip\baselineskip

\textbf{Theorem 2}
$(x^4+y^4)\frac{m^4+1}{2} = z^4+w^4$ has infinitely many parametric solutions,\\
where $m$ is arbitrary odd number.

Proof.\\

We use an identity $$p(t+1)^4+p(t)^4=(t^2+at+b)^2+(ct^2+dt+e)^2$$
where $$(a,b,c,e,p)=(d+2, \frac{1}{2d}+1, 1+d, \frac{1}{2d}, 1+d+\frac{1}{2d^2})$$\\

So, we look for the integer solutions $z^2 = t^2+(d+2)t+\frac{1}{2d}+1$ and  $w^2 = (1+d)t^2+dt+\frac{1}{2d}$.
By parameterizing the first equation and substituting the result to second equation, then we obtain quartic equation.

Let $d=m^2-1$ and $U=1/k$ then
\begin{equation}
V^2 = (4m^6+5m^2)U^4+(-8m^4-4)U^3-2m^2U^2+4U+m^2
\end{equation}

This quartic equation is birationally equivalent to an elliptic curve below.

$$E: Y^2-4YX/m+(8m+16m^5)Y = X^3-2(m^4+2)X^2/(m^2)+(-20m^4-16m^8)X+104m^6+32m^{10}+80m^2$$

Transformation is given, 
\begin{align*}
&U = \frac{2m^2X-4m^4-8}{Ym}  \\
&V = \frac{m^4X^3-6m^6X^2-12X^2m^2+28m^8X+32m^4X+32X+16Ym-16m^9Y+16m^{12}X-32m^{14}-104m^{10}-80m^6}{Y^2m^3} \\
&X = \frac{2mV+2m^2+4U}{U^2} \\
&Y = \frac{4m^3V+4m^4+8m^2U-4m^4U^2-8U^2}{U^3m}
\end{align*}

It has a point $P(X,Y)=( \frac{2(m^4+2)}{m^2}, \frac{-16(-1+m^8)}{m^3} )$.\\
According to the Nagell-Luts theorem, since the point $P$ is a point of infinite order then we can obtain infinitely many rational points on $E$.\\
Thus we can obtain infinitely many rational solutions of $(2)$ using group law. \\
Therefore $(x^4+y^4)\frac{m^4+1}{2} = z^4+w^4$ has infinitely many parametric solutions. \\
The proof is completed. \\

When $m$ is even number, problem links to $(m^4+1)(x^4+y^4) = 2(z^4+w^4)$. 
\vskip\baselineskip

Example

$n=(m^4+1)/2$ \\
$m$ is odd number.

\begin{align*}
&x = m^4+4m^2-1 \\
&y = m^4-4m^2-1 \\
&z = 3m^4+1 \\
&w = (m^4+3)m
\end{align*}
\vskip\baselineskip

\begin{align*}
&x = m^{12}+12m^{10}-19m^8+40m^6+19m^4+12m^2-1 \\
&y = m^{12}-12m^{10}-19m^8-40m^6+19m^4-12m^2-1 \\
&z = 5m^{12}-27m^8-41m^4-1 \\
&w = m(m^{12}+41m^8+27m^4-5)
\end{align*}

\newpage

\centerline{\huge4.\quad Example for $41(x^4+y^4) = z^4+w^4$}
\vskip\baselineskip
We show a numerical example of first method.
 
Let $(a,b,m)=(4,5,3)$ then $(x^4+y^4)(a^2+b^2)=z^4+w^4$ is reduced to 
$$V^2 = U^4-108U^3-18U^2+996U+409$$

Quartic equation is birationally equivalent to an elliptic curve $E$.

$$E: Y^2+XY+Y=X^3-X^2-27X+26$$

Transformation is given, 
\begin{align*}
U &=  \frac{4Y+29X-10}{X-46}\\
V &=  \frac{16X^3-2208X^2-3392X+13744-9840Y}{(X-46)^2} \\
X &=  \frac{V+U^2-54U+5}{32} \\
Y &=  \frac{UV+U^3-83U^2+99U-29V+175}{128}
\end{align*}

The rank of elliptic curve $E$ is $1$ and has generator $P(X,Y)=(6,-11)$ using SAGE. \\
We can find infinitely many rational points on the curve E using the group law.\\
Thus the multiples nP, n = 2, 3, ...give infinitely many points as follows.\\

\begin{table}[hbtp]
\centering
  \caption{Solutions of $41(x^4+y^4) = z^4+w^4$}
  \label{table:data_type}
  \begin{tabular}{cllll}
    \hline
   nP & x  & y  &  z  & w   \\
    \hline \hline

$2P$ & 29 & 11 & 63 & 61 \\
$3P$ & 17909 & 5149 & 37623 & 38699 \\
$4P$ & 229422601 & 214213319 & 663306603 & 282177719 \\
$5P$ & 81840455152441 & -86237007592439 & 252933880274523 & 61172008172039 \\
    \hline
  \end{tabular}
\end{table}

\newpage

\centerline{\huge5.\quad Example for $17(x^4+y^4) = z^4+w^4$}
\vskip\baselineskip

The previous methods can't give the rational solution of  $17(x^4+y^4) = z^4+w^4$.
Hence we show the third method of giving the rational solution to the equation $(x^4+y^4)(a^2+b^2) = z^4+w^4$.
We consider the simultaneous equation as follows.
\begin{displaymath}
\left\{
\begin{array}{l}
z^2-w^2-5x^2-3y^2=2t(4x^2+w^2-y^2) \\
(w^2+z^2+5x^2+3y^2)t=-(x^2+w^2+4y^2)
\end{array}
\right.
\end{displaymath}

Eliminating $z^2$,  we get
\begin{equation}
(1+10t+8t^2)x^2+(4+6t-2t^2)y^2+(2t^2+1+2t)w^2 = 0
\end{equation}

Parametric solution of $(3)$ is $(x,y,w)=(2(k-2)(k-3), (k-1)(3k-7), -(25+5k^2-22k)).$
\begin{align*}
z^2&=\frac{-w^2+17y^2}{4}  \\ 
   &=\frac{4(2k^2-5k+1)(4k^2-15k+13)}{(-5+k^2)^2}
\end{align*}
Hence cosider the quartic equation $(4)$

\begin{equation}
V^2=8U^4-50U^3+105U^2-80U+13
\end{equation}

Quartic equation $(4)$ has a rational point $Q(U,V)=(2,1)$, then this quartic equation is birationally equivalent to an elliptic curve below.

$E: Y^2 =X^3-91X+330$ with rank is $1$ and generator=$(7, -6)$.

Transformation is given, 
\begin{align*}
U &=  \frac{6X-36+2Y}{Y+2X-12}\\
V &=   \frac{-18X^2-114+X^3+91X}{(Y+2X-12)^2}\\
X &=  \frac{2V+6-U^2}{U^2-4U+4} \\
Y &=  \frac{12V-108+180U-90U^2-4VU+14U^3}{U^3-6U^2+12U-8}
\end{align*}

The rank of elliptic curve $E$ is $1$ and has generator $P(X,Y)=(7,-6)$ using SAGE. \\
We can find infinitely many rational points on the curve E using the group law.\\
Thus the multiples nP, n = 2, 3, ...give infinitely many points as follows.\\

\begin{table}[hbtp]
  \caption{Solutions of $17(x^4+y^4) = z^4+w^4$}
  \label{table:data_type}
  \begin{tabular}{cllll}
    \hline
   nP & x  & y  &  z  & w   \\
    \hline \hline

$2P$  & 3120 & 1921 &  2242 & 6529\\
$3P$  & 18418554 & 88538885 &  176117272 & 95896333\\
$4P$  & 87733253643360 & 108376421998081 & 198203611434238 & 206237591201281\\
$5P$  & 12509563104278834954874 & 6446124521923428875525 &  3117838409641509334568 & 25836199364300466735373\\
    \hline
  \end{tabular}
\end{table}

\newpage

\centerline{\huge6.\quad Example for $97(x^4+y^4) = z^4+w^4$}
\vskip\baselineskip

The previous methods can't give the rational solution of  $97(x^4+y^4) = z^4+w^4$.
Hence we show the fourth method of giving the rational solution to the equation $(x^4+y^4)(a^2+b^2) = z^4+w^4$.
According to Richmond's theorem\cite{d}, existence of solution for diophantine equation $ax^4 + by^4 + cz^4 + dw^4 = 0$ are known if $abcd$ is square number.
He proved other solution is derived from a known solution.
Repeting this process, we can obtain infinitely many integer solutions.
We use a known solution $(x,y,z,w)=(112,71,10,37)$ obtained by brute force as follows.\\
Let $x=pt+112, y=qt+71, z=rt+10, w=t+37$.\\
We obtain
\begin{align*}
&(q^4+p^4-97r^4-97)t^4+(-3880r^3+284q^3-14356+448p^3)t^3 \\
&+(-796758-58200r^2+30246q^2+75264p^2)t^2+(-19653364+5619712p-388000r+1431644q)t=0
\end{align*}

In order to set the coefficient of $t$ and $t^2$ to zero, we obtain\\
\begin{align*}
&q = \frac{-135744p}{42103}+\frac{491693}{42103}\\
&r = \frac{-22422}{2965}+\frac{7672p}{2965}
\end{align*}

Hence we obtain

$$t = \frac{-1238771307200}{68835707869p-174887242544}$$
\vskip\baselineskip

Take p=1 then $(q,r)=(\frac{355949}{42103},\frac{-2950}{593})$.
\vskip\baselineskip

Finally, we obtain 
$(x,y,z,w)=(174887242544, 240033770927, 68026751110, 68835707869).$\\

Thus we can obtain other integer solution by using this new solution as a known solution.\\
According to Richmond's theorem, all solutions of Table 4 have infinitely many integer solutions.

\newpage

\begin{center}
Smallest: $Min\ n(x^4 + y^4)$  \\
Brute force search range: $n<1000, (x,y,z,w)<50000$  \\
$n\equiv 1,2,8,9\mod 16$ is to be searched.
\end{center}

\begin{table}[hbtp]
\centering
  \caption{Solutions of $n(x^4+y^4) = (z^4+w^4)$}
  \label{table:data_type}
  \begin{tabular}{rrrrr}
    \hline
   n & x  & y  &  z  & w \\
    \hline \hline

2 & 7 & 20 & 21 & 19 \\
8 & 19 & 21 & 40 & 14 \\
17 & 5 & 6 & 13 & 8 \\
41 & 1 & 1 & 3 & 1 \\
82 & 219 & 320 & 1011 & 247 \\
97 & 10 & 37 & 112 & 71 \\
113 & 1 & 2 & 6 & 5 \\
136 & 1 & 1 & 4 & 2 \\
137 & 29 & 5 & 99 & 31 \\
146 & 1 & 2 & 7 & 3 \\
178 & 1 & 2 & 7 & 5 \\
193 & 18 & 43 & 159 & 80 \\
226 & 1 & 8 & 31 & 7 \\
241 & 1 & 2 & 8 & 1 \\
257 & 4 & 15 & 52 & 49 \\
313 & 1 & 1 & 5 & 1 \\
328 & 7 & 3 & 30 & 8 \\
337 & 15 & 34 & 147 & 26 \\
353 & 1 & 1 & 5 & 3 \\
386 & 1 & 2 & 9 & 1 \\
401 & 1 & 2 & 9 & 4 \\
433 & 8 & 13 & 53 & 50 \\
466 & 2 & 43 & 181 & 151 \\
482 & 5 & 6 & 31 & 7 \\
521 & 9 & 1 & 43 & 1 \\
562 & 10 & 11 & 61 & 7 \\
577 & 175 & 188 & 929 & 848 \\
578 & 24 & 37 & 187 & 85 \\
593 & 1 & 2 & 10 & 3 \\
626 & 16 & 29 & 141 & 97 \\
641 & 16 & 19 & 97 & 78 \\
673 & 161 & 26 & 820 & 139 \\
706 & 17 & 28 & 149 & 13 \\
712 & 15 & 13 & 86 & 36 \\
761 & 7 & 3 & 37 & 11 \\
776 & 1 & 1 & 6 & 4 \\
802 & 1 & 10 & 53 & 19 \\
857 & 7 & 5 & 39 & 23 \\
866 & 1 & 2 & 11 & 3 \\
881 & 7 & 31 & 169 & 9 \\
898 & 1 & 2 & 11 & 5 \\
953 & 2041 & 2021 & 12975 & 7999 \\
977 & 3 & 10 & 56 & 11 \\
                 
    \hline
  \end{tabular}
\end{table}

\newpage
\centerline{\huge7. Final remarks}
\vskip\baselineskip

There is an interesting relationship between $8(19^4+21^4)=14^4+40^4$ and $2(7^4+20^4)=19^4+21^4$ \\
Multiply both sides of  $2(7^4+20^4)=19^4+21^4$ by $2^3$, then we obtain $14^4+40^4=8(19^4+21^4)$.\\
Thus an equation $n(x^4+y^4)=z^4+w^4$ can be transform  $(nx)^4+(ny)^4 = n^3(z^4+w^4)$.\\
We could not find a solution of $n(x^4+y^4)=z^4+w^4$ for $n=34,40,50,65,73,89$ whre $n<100$.
They are not the form $(m^4+n^4)/2$.
I don't know if they don't have a solution to begin with, or if they have a large solution.
\vskip3\baselineskip


\begin{thebibliography}{99}
\bibitem{a} Ajai Choudhry, Iliya Bluskov and Alexander James, A quartic diophantine equation inspired by Brahmagupta's identity, 2020, arXiv.org
\bibitem{b}F. Izadi and K. Nabardi, Diophantine equation $X^4 + Y^4 = 2(U^4 + V^4 )$, Math. Slovaca 66 (2016),
\bibitem{c} Ali S. Janfada , Kamran Nabardi, ON DIOPHANTINE EQUATION $x^4 + y^4 = n(u^4 + v^4 )$, December 2019, Mathematica Slovaca 69(6)
\bibitem{d} H. W. Richmond, On the Diophantine equation $F  = ax^4 + by^4 + cz^4 + dw^4$, the product $abcd$ being square number, J. Lond. Math. Soc., 19 (1944) 
\end{thebibliography}
\end{document}